\RequirePackage{fix-cm}
\documentclass{svjour3}                     
\usepackage{graphicx}
 \usepackage{amsfonts}
\newtheorem{lem}{Lemma}

\newtheorem{thm}{Theorem}
\newtheorem{asp}{Assumption}

\newtheorem{mydef}{Definition}
\begin{document}

\title{A Note On the Rank of the Optimal Matrix in Symmetric Toeplitz Matrix Completion Problem
}


\author{Xihong Yan         \and
       Jiahao Guo           \and
        Yi Xu 
}


\institute{Xihong Yan \at
              Department of Mathematics, Taiyuan Normal University, Jinzhong 030619, Shanxi Province, P.R. China. \\
              \email{xihong1@e.ntu.eu.sg}           
           \and
            Jiahao Guo \at
              Department of Mathematics, Taiyuan Normal University, Jinzhong 030619, Shanxi Province, P.R. China. \\
              \email{1015369806@qq.com} \\
            Xu Yi (corresponding author)\at
              Institute of Mathematics, Southeast University, Nanjing 210096, Jiangsu Province,  P.R. China. \\
              \email{yi.xu1983@hotmail.com}
}


\maketitle

\begin{abstract}
We consider the symmetric Toeplitz matrix completion problem, whose matrix under consideration possesses specific row and column structures. This problem, which has wide application in diverse areas, is well-known to be computationally NP-hard. This note provides an upper bound on the objective of minimizing the rank of the symmetric Toeplitz matrix in the completion problem based on the theorems from the trigonometric moment problem and semi-infinite problem. We prove that this upper bound is less than twice the number of linear
constraints of the Toeplitz matrix completion problem. Compared with previous work in the literature, ours is one of the first efforts to investigate the bound of the objective value of the Toeplitz matrix completion problem.
\end{abstract}

\section{Introduction}
The origins of the well-known low rank matrix completion problem can be traced back to the works of Prony in 1795. It has gained tremendous popularity recently due to wide applications in various fields including machine learning \cite{Argyriou2007}, compressed sensing\cite{Candes2009,Chen2013}, system identification and control \cite{Fazel2013,Liu2010}, computer vision \cite{Tomasi1992}, and so on. A lot of work has been developed to address methods to solve the matrix completion problem, such as the singular value thresholding algorithm \cite{Cai2010}, the accelerated proximal gradient algorithm \cite{Toh2010}, and the matrix factorization based approach \cite{Sun2016}. For many other methods are well documented in \cite{Fazel2002,Hu2012,Goldfarb2011}.

In many engineering and statistical applications, the matrices under consideration are often structured. For example, in statistical signal processing, the covariance matrix of a
stationary random process usually has a Toeplitz structure. In fact, due to the important role of a Toeplitz matrix on real-world
problems arising in signal and image sciences, numerous research has been conducted to discuss the Toeplitz matrix \cite{Shaw1998,Wang2017}. As a consequence, the Toeplitz matrix completion problem, recovering an unknown low-rank or approximately low-rank Toeplitz matrix from a sampling of its entries, has become an extremely important issue. The problem can be characterized mathematically as the following linearly constrained minimization problem,
\begin{eqnarray}\label{mb1}
  \min &r(X) &\nonumber\\
   s.t.&X_{ij} =M_{ij}, &\forall (i,j)\in\Omega,\\
 &X \mbox{\ is\  Toeplitz}.& \nonumber
\end{eqnarray}
where $r(X)$ is the rank of $X$, $M\in \mathbb R^{n_1\times n_2}$ is a given Toeplitz matrix for which only a subset of its entries
$M_{ij},(i,j)\in \Omega\subset\{1,2,...n_1\}\times\{1,2,...n_2\}$ are specified.
However, methods for the Toeplitz matrix completion problem (\ref{mb1}) has not been adequately addressed in the previous studies. Several evolutionary algorithms were proposed by Wang based on the nuclear norm model \cite{Wang2015,Wang2016}. Wang first introduced the augmented lagrange multiplier algorithm \cite{Wang2016} and the singular value thresholding algorithm \cite{Wang2015} for solving the Toeplitz matrix completion problem, where each iteration kept a feasible Toeplitz structure. The major drawbacks for the existing methods for the Toeplitz matrix completion problem are that there is not a theoretical guarantee to correctly recover the underlying Toeplitz matrix due to the
non-convexity of the problem and they are highly sensitive to the choice of some parameters which are dependent on a guess of the rank of the recovered matrix. The range of the rank of the optimal Toeplitz matrix plays very important roles in terms of improving performance of algorithms for solving problem (\ref{mb1}). Therefore, the purpose of this paper is to provide an upper bound on the rank of the optimal matrix for problem (\ref{mb1}).

Most literature about estimation of the rank of a matrix concentrates on the rank of the optimal matrix of a positive semidefinite program \cite{Fawzi2015,Gouveia2015}. Compared with previous work in the literature, ours is one of the first efforts to address an upper bound on the rank of the recovered matrix in the symmetric Toeplitz matrix completion problem, where the procedure is as follows: First, we describe the feasible region of problem (\ref{mb1}) as a positive semidefinite system by using simple matrix decomposition; Then, we adopt a theorem from the trigonometric moment problem to convert the positive semidefinite system to a semi-infinite problem; Finally, based on the optimality conditions for the semi-infinite problem, we establish the theorem of an upper bound for the minimum objective of problem (\ref{mb1}).

The rest of this paper is organized as follows. In Section 2, we define our notation and
give some preliminaries for the subsequent analysis. Section 3 presents
our main results on an upper bound on the rank of the optimal matrix in the symmetric Toeplitz matrix completion problem.
Finally, we give some concluding remarks in Section 4.

\section{Preliminaries}
In this section, some notations and basic definitions are summarized, which will be used in the remaining part of the paper. We then give a well-known result on the solution of the trigonometric moment problem that will play central roles in the later analysis.

\textbf{Notation.} Let $\mathbb R^{n}$ be an $n$-dimensional Euclidean space, $\mathbb R^{n\times n}$ be the set of $n\times n$ real matrices. $r(T)$ denotes the rank of a matrix $T$. A matrix $T$ is symmetric positive definite (resp. positive semidefinite) and is denoted by $T\succ0$ (resp. $T\succeq 0$). $i$ represents the unit imaginary number. $|A|$ is the number of elements of set $A$. The symbol $^T$ represents the transpose. For a measure $u(t)$, its supporting set is defined as $supp(\mu(t)):=\{t\in\mathcal{T} |\mu(t)>0\}$.

\begin{mydef}\label{def-1}
 An $n\times n$ real Toeplitz matrix is the name for a matrix with the following shape,
   \begin{eqnarray*}{
\left(\begin{array}{ccccc}
  x_{0}& x_{1} &\cdots & x_{n-2}& x_{n-1} \\
  x_{-1}  & x_{0} & \cdots & x_{n-3} & x_{n-2} \\
  \vdots & \vdots & \ddots& \vdots& \vdots \\
  x_{-n+2}& x_{-n+3} & \cdots & x_{0} & x_{1}\\
   x_{-n+1}& x_{-n+2} & \cdots & x_{-1} & x_{0}
\end{array}
\right).}
\end{eqnarray*}
\end{mydef}
Noted that an $n \times n$ Toeplitz matrix is an $n \times n$ matrix whose
entries $x_{kj}$ satisfy $x_{kj}=x_{j-k}$ for all $k$ and $j$, which implies that it is determined by $2n-1$ entries which are in the first row and first column. If we let $\bar{x}:=(x_{-n+1},\cdots,x_{-1},x_0,x_1,\cdots, x_{n-1})^T\in \mathbb R^{2n-1}$, it is convenient to denote the associated Toeplitz matrix as $\bar{T}(\bar{x}).$

When a Toeplitz matrix is symmetric, we obtain a symmetric Toeplitz matrix defined as follows.
\begin{mydef}\label{def-2}
A matrix with the following shape is called as an $n\times n$ symmetric Toeplitz matrix,
   \begin{eqnarray*}{
\left(\begin{array}{ccccc}
  x_{0}& x_{1} &\cdots & x_{n-2}& x_{n-1} \\
  x_{1}  & x_{0} & \cdots & x_{n-3} & x_{n-2} \\
  \vdots & \vdots & \ddots& \vdots& \vdots \\
  x_{n-2}& x_{n-3} & \cdots & x_{0} & x_{1}\\
   x_{n-1}& x_{n-2} & \cdots & x_{1} & x_{0}
\end{array}
\right).}
\end{eqnarray*}
\end{mydef}
It is clear that an $n \times n$ symmetric Toeplitz matrix is determined by $n$ entries of the first row denoted by a vector $x:=(x_{0}, x_{1} ,\cdots, x_{n-1})^T\in \mathbb R^{n}$. Hence, the corresponding symmetric Toeplitz matrix can be denoted by $T(x)$. Notice that when $X$ is a symmetric Toeplitz matrix $T(x)$, model (\ref{mb1}) can be revised as follows,
\begin{eqnarray}\label{mb2}
  \min &r(T(x)) &\nonumber\\
   s.t.&(T(x))_{kj}=M_{kj}, &\forall (k,j)\in\Omega.
\end{eqnarray}
The constraints $(T(x))_{kj}=M_{kj}, \forall (k,j)\in\Omega$ of the above model express the condition that
the recovered matrix is consistent with the observed data. These linear equality constraints can be taken the following form,
$$Bx=d,$$
where $B=(b_{k,j})_{k=1,j=0}^{m, n-1}\in \mathbb R^{m\times n}$ and $d=(d_{1}, d_{2} ,\cdots, d_{m} )^T\in \mathbb R^{m}.$ When we introduce the form of $Bx=d$ for all equality constraints in (\ref{mb2}), the resulting more general formulation is as follows,
\begin{eqnarray}\label{mb3}
 & \min &r(T(x)) \nonumber\\
   &s.t. &Bx=d.
\end{eqnarray}
Program (\ref{mb3}) is a specially structured matrix completion optimization problem since the constraint that the concerned matrix is a symmetric Toeplitz matrix is required. This requirement however comes at a price, with theoretical and practical difficulties involved in the process of solving such a problem. Hence, we present a lemma about the solution of the classical trigonometric moment problem, which is a frequently used and powerful idea for simplifying the high technical requirement.  For clearness, we first review the trigonometric moment problem. The problem is that asking whether a prescribed finite sequence $\{\alpha_0,\ldots,\alpha_{n-1}\}$ can be represented as the sequence of successive
moments of some positive measure. Specifically, for given sequence $\{\alpha_0,\ldots,\alpha_{n-1}\}$, whether there exists a positive Borel measure $\mu(t), t\in[0,2\pi]$ such that
$$\alpha_j=\frac{1}{2\pi}\int_{0}^{2\pi}e^{-ijt}\mu(dt), \ j = -n+1,\cdots,n-1,$$ where $\alpha_{-j}$ and $\alpha_{j}$ are conjugate to each other. The question is answered by previous work \cite{trm2,trm1} whose result is stated in the following lemma.

\begin{lem}\label{lem-0}
The trigonometric moment problem has a solution, \textit{i.e.}, $\{\alpha_j\}_{j=0}^{n-1}$ is a valid sequence of Fourier coefficients of some positive Borel measure $\mu(t)$ defined on $[0, 2\pi]$, if and only if there exists a Hermite Toeplitz matrix $T$ that is positive semidefinite, where
$$T=\left(
   \begin{array}{cccc}
     \alpha_0 & \alpha_1 & \cdots & \alpha_{n-1} \\
     \alpha_{-1} & \alpha_0 & \cdots & \alpha_{n-2} \\
     \vdots & \vdots & \ddots & \vdots \\
  \alpha_{-n+1} & \alpha_{-n+2} & \cdots & \alpha_0 \\
   \end{array}
 \right).$$
 \end{lem}
\section{Upper Bound for the Rank of the Optimal Matrix in Symmetric Toeplitz Matrix Completion Problem}

In this section, we focus on analysis of an upper bound of the rank of the solution of the symmetric Toeplitz matrix completion problem (\ref{mb3}). To achieve the goal, we look for a feasible point $x$ of program (\ref{mb3}) and estimate the rank of the associated Toeplitz matrix $T(x)$, which suffices to find a feasible solution of the following system and estimate the rank of its associated Toeplitz matrices,
\begin{eqnarray}\label{eq-2}
& Bw_1-Bw_2=d,\nonumber \\
  & T(w_1)\succeq 0, T(w_2)\succeq 0,
\end{eqnarray}
where $w_1=((w_1)_{0},...,(w_1)_{n-1})^T\in \mathbb R^{n}$ and $w_2=((w_2)_{0},...,(w_2)_{n-1})^T\in \mathbb R^{n}$. It is obvious that we can find a feasible point $x$ of program (\ref{mb3}) by finding a solution $(w_1,w_2)$ of system (\ref{eq-2}) and setting $x=w_1-w_2$. Moreover, we can obtain the estimation of $r(T(x))$ by investigating the rank of matrices $T(w_1)$ and $T(w_2).$

We enlarge matrix $B$ to $\bar{B}=(\bar{b}_{k,j})\in \mathbb R^{m\times(2n-1)}$ and vector $w_k$ to $\bar{w}_k (k=1,2)$ in the following way,

\begin{eqnarray}\label{b}
\bar{b}_{k,j}=\left\{
        \begin{array}{ccc}
          \frac{b_{k,-j}}{2}, &  j<0,&\\
          b_{k,0}, &  j=0,&k=1,2,...,m,\ j=-n+1,\cdots,n-1,\\
          \frac{b_{k,j}}{2}&  j>0,&\\
        \end{array}
      \right.
\end{eqnarray}
and
\begin{eqnarray}\label{w}(\bar{w}_{k})_{j}=\left\{
        \begin{array}{ccc}
          \frac{(w_k)_{-j}}{2}, &  j<0,&\\
          (w_{k})_{0}, &  j=0,&k=1,2,\ j=-n+1,\cdots,n-1.\\
          \frac{(w_{k})_{j}}{2},&  j>0,&\\
        \end{array}
      \right.
\end{eqnarray}
Then, we have the following equivalent expression of system (\ref{eq-2}),
\begin{eqnarray}\label{eq-3}
 & \bar{B}\bar{w}_1-\bar{B}\bar{w}_2=d,\nonumber\\
   & (\bar{w}_1)_j=(\bar{w}_1)_{-j},&j=1,\cdots,n-1,\\
   & (\bar{w}_2)_j=(\bar{w}_2)_{-j},&k=1,\cdots,n-1,\nonumber\\
    & \bar{T}(\bar{w}_1)\succeq 0, \bar{T}(\bar{w_2})\succeq 0.& \nonumber
\end{eqnarray}
Applying the statement of Lemma \ref{lem-0} to $\bar{w}_1$ and $\bar{w}_2$ in the above expression, that is, there exist two positive Borel-signed measures defined on $[0, 2\pi]$ such that
$$(w_1)_j=\displaystyle\frac{1}{2\pi}\int_{0}^{2\pi}e^{-ijt}\mu_1(dt), \ j=-n+1,\cdots,n-1,$$ and $$(w_2)_j=\displaystyle\frac{1}{2\pi}\int_{0}^{2\pi}e^{-ijt}\mu_2(dt),\  j=-n+1,\cdots,n-1,$$
we obtain the following system,
\begin{eqnarray}\label{Programut}
& \displaystyle\sum_{j=-n+1}^{n-1}\bar{b}_{k,j}\frac{1}{2\pi}\int_{0}^{2\pi}e^{-ijt}\mu_1(dt)-\displaystyle\sum_{j=-n+1}^{n-1}\bar{b}_{k,j}\frac{1}{2\pi}\int_{0}^{2\pi}e^{-ijt}\mu_2(dt)=d_k,\nonumber\\
&~~~~~~~~~~~~~~~~~~~~~~~~~~~~~~~~~~~~~~~~~~~~~~~~~~~~~~k=1,2,\cdots,m,\nonumber\\
 &\displaystyle\int_{0}^{2\pi}e^{-ijt}\mu_1(dt)= \int_{0}^{2\pi}e^{ijt}\mu_1(dt),\ \  j=-n+1,...,n-1,\\
 & \displaystyle\int_{0}^{2\pi}e^{-ijt}\mu_2(dt)= \int_{0}^{2\pi}e^{ijt}\mu_2(dt),\ \  j=-n+1,...,n-1,\nonumber\\
       & \mu_1(t), \mu_2(t)\in P([0, 2\pi]),&\nonumber
\end{eqnarray}
where $P([0, 2\pi])$  is a set of the finite positive Borel-signed measures defined on $[0, 2\pi].$
Simplifying program (\ref{Programut}), its equivalent expression is as follows,
\begin{eqnarray}\label{eq-41}
& \displaystyle\sum_{j=-n+1}^{n-1}\bar{b}_{k,j}\frac{1}{2\pi}\int_{0}^{2\pi}\cos(jt)\mu_1(dt)-\displaystyle\sum_{j=-n+1}^{n-1}\bar{b}_{k,j}\frac{1}{2\pi}\int_{0}^{2\pi}\cos(jt)\mu_2(dt)=d_k,\nonumber\\
&~~~~~~~~~~~~~~~~~~~~~~~~~~~~~~~~~~~~~~~~~~~~~~~~~~~~~~k=1,2,\cdots,m,\nonumber\\
   &\displaystyle\int_{0}^{2\pi}\sin(jt)\mu_1(dt)=0,\ \  j=-n+1,...,n-1,\\
 & \displaystyle\int_{0}^{2\pi}\sin(jt)\mu_2(dt)=0,\ \ j=-n+1,...,n-1,\nonumber\\
       & \mu_1(t), \mu_2(t)\in P([0, 2\pi]).&\nonumber
\end{eqnarray}
It would be a computationally onerous task to directly determine a solution of system (\ref{eq-41}). Hence, we provide an easier way to work it out. The procedure is as follows: First, we solve program (\ref{eq-4}) 

\begin{eqnarray}\label{eq-4}
  \min& \displaystyle\sum_{j=-n+1}^{n-1}(\bar{c}_1)_{j}\frac{1}{\pi}\int_{0}^{\pi}\cos(jt)\mu_1(dt)-\displaystyle\sum_{j=-n+1}^{n-1}(\bar{c}_2)_{j}\frac{1}{\pi}\int_{0}^{\pi}\cos(jt)\mu_2(dt) &\nonumber\\
  s.t. & \displaystyle\sum_{j=-n+1}^{n-1}\bar{b}_{k,j}\frac{1}{\pi}\int_{0}^{\pi}\cos(jt)\mu_1(dt)-\displaystyle\sum_{j=-n+1}^{n-1}\bar{b}_{k,j}\frac{1}{\pi}\int_{0}^{\pi}\cos(jt)\mu_2(dt)=d_k,\nonumber\\
  &~~~~~~~~~~~~~~~~~~~~~~~~~~~~~~~~~~~~~~~~~~~~~~~~~~~~~~k=1,2,\cdots,m,\nonumber\\
       & \mu_1(t), \mu_2(t)\in P([0, \pi]),&
\end{eqnarray}
to obtain its solution $(\mu_1(t),\mu_2(t))$, where $$\bar{c}_1=((\bar{c}_1)_{-n+1},...,(\bar{c}_1)_{0},...,(\bar{c}_1)_{n-1})^T\in \mathbb R^{2n-1}$$ and $$\bar{c}_2=((\bar{c}_2)_{-n+1},...,(\bar{c}_2)_{0},...,(\bar{c}_2)_{n-1})^T\in \mathbb R^{2n-1}$$ satisfy that $(\bar{c}_k)_{j}=(\bar{c}_k)_{j}, k=1,2,\ j=-n+1,...,n-1;$
Then,
$\mu_1(t)$ and $\mu_2(t)$ can be extended to measures on $[0,2\pi]$ by defining
$$ \bar{\mu}_k(t)=\left\{
  \begin{array}{ccc}
    \mu_k(t), & t\in [0,\pi],&\\
              &             &k=1,2.\\
    \mu_k(2\pi-t), &t\in [\pi, 2\pi],&\\
  \end{array}
\right.$$

The dual program of program (\ref{eq-4}) is as follows,
\begin{eqnarray}\label{eq-6}
  \min & d^Ty \nonumber\\
  s.t.& \displaystyle \bar{a}^l_0(t)\leq \sum_{k=1}^m \bar{a}_k(t)y_k\leq \bar{a}^u_0(t), &\forall t\in [0, \pi],
\end{eqnarray}
where $$ \bar{a}^l_0(t)=\displaystyle\sum_{j=-n+1}^{n-1}(\bar{c}_1)_{j}\cos(jt),$$ $$\bar{a}^u_0(t)=\displaystyle\sum_{j=-n+1}^{n-1}(\bar{c}_2)_{j}\cos(jt),$$ and $$\bar{a}_k(t)=\displaystyle\sum_{j=-n+1}^{n-1}\bar{b}_{k,j}\cos(jt), k=1,2,\cdots,m.$$

Now, we show the boundness of the feasible set of program (\ref{eq-6}) in the following lemma.

\begin{lem}\label{lem-2}
If $B$ is full rank, the feasible set of program (\ref{eq-6}) is bounded.
\end{lem}

\proof First, we prove that for $m$ distinct points $t_j\in [0, \pi](j=1,2,...m)$, $\big\{(\bar{a}_1(t_j),\bar{a}_2(t_j),\cdots,\bar{a}_m(t_j))^{T}\big\}_{j=1}^m$ is linear independent.

Based on the definition of $\bar{a}_k(t)$, it follows that
$$\big(\bar{a}_1(t), \bar{a}_2(t),...,\bar{a}_m(t)\big)=BPz,$$
where $P$ is an $n\times n$ Chebyshev
coefficient matrix satisfying $$\displaystyle\sum_{j=0}^{n-1}P_{kj}\cos^j(t)=\cos(kt),$$ and $$z=(1, \cos(t),\cos^2(t),\ldots,\cos^n(t))^T.$$
Thus,
 $$(\bar{a}_k(t_j))_{k,j=1}^m=BPZ,$$
where $Z=\left(
           \begin{array}{cccc}
             1 & 1 & \cdots & 1 \\
             z_1 & z_2 & \cdots & z_m \\
             z_1^2 & z_2^2 & \cdots & z_m^2 \\
             \vdots & \vdots & \vdots & \vdots \\
             z_1^m&z_2^m &\cdots &z_m^m\\
           \end{array}
         \right)
$ and $z_{j}=\cos(t_j), j=1,2,...,m.$

Since $\{t_j\in [0, \pi]\}_{j=1}^m$ are distinct, $\{z_j\}_{j=1}^m$ are distinct. Therefore, $Z$ is nonsingular based on the property of Vandermonde matrix, which together with the facts that $P$ and $B$ are full rank, shows that $BPZ$ is full rank. This says there exist $m$ points $t_j\in [0, \pi]$ such that $\big\{(\bar{a}_1(t_j),\bar{a}_2(t_j),\cdots,\bar{a}_m(t_j))^{T}\big\}_{j=1}^m$ is linear independent.

Then, we show that the feasible set of program (\ref{eq-6}) is bounded. If the feasible set of program (\ref{eq-6}) denoted $C$ is unbounded, there is a direction $h\neq0$ such that for every $\lambda$ and $y\in C$, $y+\lambda h\in C.$ This says, for any $\lambda$, we have
$$\bar{a}^l_0(t_{j})\leq \sum_{k=1}^m \bar{a}_k(t_{j})(y_k+\lambda h_{k})\leq \bar{a}^u_0(t_{j}),\ \  j=1,2,...,m,$$
which implies that
\begin{eqnarray}\label{eq-fc}
\sum_{k=1}^m \bar{a}_k(t_{j})h_{k}=0,&  j=1,2,...,m.
\end{eqnarray}
Recall the linear independence of $\big\{(\bar{a}_1(t_j),\bar{a}_2(t_j),\cdots,\bar{a}_m(t_j))^{T}\big\}_{j=1}^m$, we have that system (\ref{eq-fc}) only has zero solution, which leads to $h=0$. It contradicts with the fact that $h\neq 0.$ Hence, the conclusion of this lemma is held. $\qed$

It is well known that if program (\ref{eq-6}) satisfies the slater condition and has finite optimal objective value, the objective values of program (\ref{eq-6}) and program (\ref{eq-4}) are equal. Hence, here we are making a reasonable assumption to make the slater condition hold for program (\ref{eq-6}).
\begin{asp}\label{ass-1}
$\bar{a}^u_0(t)>0$ and $\bar{a}^l_0(t)<0, \forall t\in [0,\pi].$
\end{asp}
Since $\bar{c}_1$ and $\bar{c}_2$ can be set an arbitrary vector with the element indexed by $j$ being the same as the element indexed by $-j$, for example, $(\bar{c}_1)_0=1$, $(\bar{c}_2)_0=-1$, and $(\bar{c}_k)_j=0 \ (k=1,2,\ j\neq0)$, this assumption is valid for program (\ref{eq-6}). Under this assumption, $0$ is an interior feasible point of program (\ref{eq-6}), that is, the slater condition hold for program (\ref{eq-6}).

With the above pavements, we are now ready to state the main result about an upper bound on the rank of the optimal matrix in the symmetric Toeplitz matrix completion problem.

\begin{thm}\label{thm-1}
 Suppose Assumption \ref{ass-1} holds. If $B$ is full rank and the number of its rows is $m$, then the objective value of program (\ref{mb3}) is smaller than $2m$, that is, $v(\ref{mb3})\leq 2m.$
\end{thm}

\proof Under the Assumption \ref{ass-1} and the conclusion of Lemma \ref{lem-2}, the solution set of program (\ref{eq-6}) is not empty and finite. Let $y^*$ be an optimal solution of program (\ref{eq-6}).
 Then there exists a lagrange multiplier $(\mu_1^*(t),\mu^*_2(t)) (t\in P([0,\pi]))$ satisfying the following complementary slackness conditions,
 $$\int_{0}^{\pi}\big(\bar{a}^u_0(t)-\sum_{k=1}^m \bar{a}_k(t)y^*_k\big)\mu^*_1(dt)=0,
 $$ and $$\int_{0}^{\pi}\big(\sum_{k=1}^m \bar{a}_k(t)y^*_k- \bar{a}^l_0(t)\big)\mu_2^*(dt)=0.$$
This says if $|supp(\mu_1^*(t))|+|supp(\mu_2^*(t))|>m$, there exist $l (l=l_1+l_2>m)$ points $t_1,t_2,...t_{l_{1}},...t_{l}$ such that $$\bar{a}^u_0(t_j)=\sum_{k=1}^m \bar{a}_k(t_j)y^*_k,\ \  j=1,2,\cdots,l_1,$$ and
$$\bar{a}^l_0(t_j)=\sum_{k=1}^m \bar{a}_k(t_j)y^*_k,\ \  j=l_1,l_1+1,\cdots,l.$$
Then, we can select subsets $\Delta_1$ from $\{1,\cdots, l_1\}$ and $\Delta_2$  from $\{l_1,l_1+1,\cdots,l\}$ such that $|\Delta_1|+|\Delta_2|\leq m$,
 \begin{eqnarray}\label{eq-t2}\bar{a}^u_0(t_j)=\sum_{k=1}^m \bar{a}_k(t_j)y^*_k, &j\in \Delta_1,
 \end{eqnarray}
\begin{eqnarray}\label{eq-t3}\bar{a}^l_0(t_j)=\sum_{k=1}^m \bar{a}_k(t_j)y^*_k, &j\in \Delta_2,
\end{eqnarray}
and $\big\{(\bar{a}_k(t_j)\}_{(k=1,...,m,\ j\in\Delta_1\cup \Delta_2)}$ is linear independent.
It is obvious that equations (\ref{eq-t2}) and (\ref{eq-t3}) with linear independence of $\big\{(\bar{a}_k(t_j)\}_{(k=1,...,m,\ j\in\Delta_1\cup \Delta_2)}$ hold for the case where $|supp(\mu_1^*(t))|+|supp(\mu_2^*(t))|\leq m$. In a word, we meet these requirements for program (\ref{eq-6}) and define two measures as follows,
$$ \mu_1^{**}(t)=\left\{
  \begin{array}{cc}
    \mu_1^*(t), & t=t_j, j\in  \Delta_1,\\
    0, & else, \\
  \end{array}
\right.$$
and
$$ \mu_2^{**}(t)=\left\{
  \begin{array}{cc}
    \mu_2^*(t), & t=t_j, j\in  \Delta_2,\\
    0, & else, \\
  \end{array}
\right..$$
According to the dual theorem, we obtain that $(\mu_1^{**}(t),\mu_2^{**}(t))$ is a solution of program (\ref{eq-4}).
By symmetrically extending $\mu_1^{**}(t)$ and $\mu_2^{**}(t)$ to measures on $[0,2\pi]$, respectively, we have a solution of system (\ref{eq-41}) as follows,
$$ \bar{\mu}_1(t)=\left\{
  \begin{array}{cc}
    \mu^{**}_1(t), & t\in [0, \pi],\\
    \mu^{**}_1(2\pi-t), & t\in[\pi, 2\pi], \\
  \end{array}
\right.$$
and
$$ \bar{\mu}_2(t)=\left\{
  \begin{array}{cc}
    \mu^{**}_2(t), & t\in [0, \pi],\\
    \mu^{**}_2(2\pi-t), & t\in[\pi, 2\pi]. \\
  \end{array}
\right.$$

Let $$(\bar{w}_1)_j=\frac{1}{2\pi}\int_0^{2\pi}\cos(jt)\bar{\mu}_1(dt),\ \ j=-n+1,\cdots,n-1,$$ and $$(\bar{w}_2)_j=\frac{1}{2\pi}\int_0^{2\pi}\cos(jt)\bar{\mu}_2(dt), \ \ j=-n+1,\cdots,n-1.$$
We plug $\bar{w}_1$ and $\bar{w}_2$ into system (\ref{Programut}) and obtain that $\bar{B}\bar{w}_1-\bar{B}\bar{w}_2=d$. Recalling the relationship (\ref{b}) between $B$ and $\bar{B}$ and (\ref{w}) between $w_{k}$ and $\bar{w}_{k} (k=1,2)$, we have $Bw_1-Bw_2=d$ and $r(\bar{T}(\bar{w}_{k}))=r(T(w_{k})).$ Let $x=w_1-w_2$, then $x$ is a feasible solution of program (\ref{mb3}) since $B(w_1-w_2)=d.$

By computing the above integrals, we have
 \begin{eqnarray}\label{int1}(\bar{w}_1)_j=\left\{
                                \begin{array}{cc}
                                  \displaystyle\sum_{t_k\in \Delta_1, t_k\neq \pi} 2u^*_1(t_k)\cos(jt_k)+u^*_1(\pi)\cos(j\pi), &  \mbox{if~~} \pi \in\{t_k, k\in \Delta_1\},\\
                                  \displaystyle\sum_{t_k\in \Delta_1} 2u^*_1(t_k)\cos(jt_k), &  else, \\
                                \end{array}
                              \right.
\end{eqnarray}
and
 \begin{eqnarray}\label{int2}(\bar{w}_2)_j=\left\{
                                \begin{array}{cc}
                                  \displaystyle\sum_{t_k\in \Delta_2, t_k\neq \pi} 2u^*_2(t_k)\cos(jt_k)+u^*_2(\pi)\cos(j\pi), &  \mbox{if~~} \pi \in\{t_k, k\in \Delta_2\},\\
                                  \displaystyle\sum_{t_k\in \Delta_2} 2u^*_2(t_k)\cos(jt_k), &  else. \\
                                \end{array}
                              \right.
\end{eqnarray}

Now, we explore an upper bound for the rank of Toeplitz matrix $T(x)$ by investigating the ranks of Toeplitz matrices $\bar{T}(\bar{w}_1)$ and $\bar{T}(\bar{w}_2)$. For the ease of exposition, we use $\{t_{k_{1}}, t_{k_2},...,t_{k_{|\Delta_1|}}\}$ (resp. $\{t_{s_{1}}, t_{s_{2}},...,t_{s_{|\Delta_2|}}\}$) to represent set $\{t_k| k\in \Delta_1\}$ (resp. $\{t_k| k\in \Delta_2\}$). Using Definition \ref{def-1}, decomposition of a symmetric Toeplitz matrix, and expressions (\ref{int1}) and (\ref{int2}), we have $\bar{T}(\bar{w}_1)$ and $\bar{T}(\bar{w}_2)$ as follows,

$$\bar{T}(\bar{w}_1)=V_1D_1V_1^*\  \mbox{and}\  \bar{T}(\bar{w}_2)=V_2D_2V_2^*,$$
where

$$V_1=\left(
             \begin{array}{ccccccc}
               1 & 1 & \cdots & 1 &1&\cdots&1\\
               e^{it_{k_{1}}} & e^{it_{k_{2}} }& \cdots & e^{it_{k_{|\Delta_1|}}}& e^{i(2\pi-t_{k_{|\Delta_1|}})} & \cdots &e^{i(2\pi-t_{k_{1}})} \\
               e^{i2t_{k_{1}}} & e^{i2t_{k_2}} & \cdots & e^{i2t_{k_{|\Delta_1|}}}& e^{i2(2\pi-t_{k_{|\Delta_1|}})}  & \cdots &e^{i2(2\pi-t_{k_{1}})}\\
               \cdots & \cdots & \cdots & \cdots&\cdots & \cdots & \cdots  \\
                e^{i(n-1)t_{k_{1}}} & e^{i(n-1)t_{k_2}} & \cdots & e^{i(n-1)t_{k_{|\Delta_1|}}}& e^{i(n-1)(2\pi-t_{k_{|\Delta_1|}})}  & \cdots &e^{i(n-1)(2\pi-t_{k_{1}})}\\
             \end{array}
           \right),$$
$$D_1=\left(
             \begin{array}{cccccccc}
               u^*_1(t_{k_1}) &  &  & & &&\\
                & u^*_1(t_{k_2})&  & & &&\\
               &  & \ddots & &&&\\
             & &  &  u^*_1(t_{k_{|\Delta_1|}})&&&\\
              & &  & & u^*_1(t_{k_{|\Delta_1|}})&&&\\
              & &  & & &\ddots&&\\
              & &  & & &&u^*_1(t_{k_1})&\\
             \end{array}
           \right),
$$

$$V_2=\left(
             \begin{array}{ccccccc}
               1 & 1 & \cdots & 1 &1&\cdots&1\\
               e^{it_{s_{1}}} & e^{it_{s_{2}} }& \cdots & e^{it_{s_{|\Delta_2|}}}& e^{i(2\pi-t_{s_{|\Delta_2|}})} & \cdots &e^{i(2\pi-t_{s_{1}})} \\
               e^{i2t_{s_{1}}} & e^{i2t_{s_2}} & \cdots & e^{i2t_{s_{|\Delta_2|}}}& e^{i2(2\pi-t_{s_{|\Delta_2|}})}  & \cdots &e^{i2(2\pi-t_{s_{1}})}\\
               \cdots & \cdots & \cdots & \cdots&\cdots & \cdots & \cdots  \\
                e^{i(n-1)t_{s_{1}}} & e^{i(n-1)t_{s_2}} & \cdots & e^{i(n-1)t_{s_{|\Delta_2|}}}& e^{i(n-1)(2\pi-t_{s_{|\Delta_2|}})}  & \cdots &e^{i(n-1)(2\pi-t_{s_{1}})}\\
             \end{array}
           \right),$$
$$D_2=\left(
             \begin{array}{cccccccc}
               u^*_2(t_{s_1}) &  &  & & &&\\
                & u^*_2(t_{s_2})&  & & &&\\
               &  & \ddots & &&&\\
             & &  &  u^*_2(t_{s_{|\Delta_2|}})&&&\\
              & &  & & u^*_2(t_{s_{|\Delta_2|}})&&&\\
              & &  & & &\ddots&&\\
              & &  & & &&u^*_2(t_{s_1})&\\
             \end{array}
           \right),
$$
and $V_{j}^*$ and $V_{j}(j=1,2)$ are conjugate transposes.
We can readily observe $r(\bar{T}(\bar{w}_1))\leq 2|\Delta_1|$ and $r(\bar{T}(\bar{w}_2))\leq 2|\Delta_2|$. Then, we have
\begin{eqnarray}
  r(T(x))&&=r(T(w_1-w_2))\leq r(T(w_1)+r(T(w_2)=r(\bar{T}(\bar{w}_1))+r(\bar{T}(\bar{w}_2))\nonumber\\
&&\leq 2|\Delta_1|+2|\Delta_2|\leq 2l_1+2l_2\leq 2m,\nonumber
\end{eqnarray}
which immediately shows the conclusion of this theorem, that is, $v(\ref{mb3})\leq 2m.\qed$
\section{Conclusion}
The main purpose of this paper is to derive a theoretical result on an upper bound on the rank of the symmetric Toeplitz matrix in the completion problem based on the theorems from the trigonometric moment problem and semi-infinite problem. We prove that this upper bound is only dependent of the number of linear
constraints of the Toeplitz matrix completion problem, that is, it is less than twice the number of linear
constraints of the problem. This result reduces the noteworthy requirement of guessing the rank of the objective matrix in existing methods for solving such problem.

Our paper sheds light on the upper bound of the rank of the concerned Toeplitz matrix. We realize that a
practicable range for the rank is complicated to derive in our current context. Actually, we test thousands of numerical examples, where $B$ and $d$ are randomly generated, to try to find a low bound of the rank in the sense of average. We find that the low bound is the number of linear constraints for all test examples. However, we have not obtain any theoretical result. Hence, finding a new mechanism to provide a low bound for the rank of the concerned matrix in completion problems is a future work direction.

\bigskip\medskip
\noindent{\bf  Acknowledgements:} \small{ X. Yan was supported in part by the STIP of Higher Education Institutions in Shanxi No. 201802103 and NSF of Shanxi province No. 201801D121022. Y. Xu was supported in part by National Natural Science Foundation of China No. 11501100, 11571178, 11671082 and 11871149.

\end{document}